\documentclass{amsart}

\usepackage{amssymb}
\usepackage{graphics}
\usepackage[all]{xy}
\CompileMatrices

\newtheorem{introthm}{Theorem}

\newtheorem{theorem}{Theorem}[section]
\newtheorem{lemma}[theorem]{Lemma}
\newtheorem{proposition}[theorem]{Proposition}
\newtheorem{corollary}[theorem]{Corollary}

\theoremstyle{definition}

\newtheorem{example}[theorem]{Example}

\def\div{{\rm div}}

\def\quot{/\!\!/}
\def\mal{\! \cdot \!}

\def\tq{\mathbin{\!{/ \! \! \lower 3pt \hbox{$\scriptscriptstyle{\rm
          tq}$}}\!}}

\def\rq#1{\widehat{#1}}
\def\t#1{\widetilde{#1}}
\def\b#1{\overline{#1}}

\def\CC{{\mathbb C}}
\def\KK{{\mathbb K}}
\def\TT{{\mathbb T}}
\def\ZZ{{\mathbb Z}}
\def\RR{{\mathbb R}}

\def\QQ{{\mathbb Q}}

\def\CDiv{{\rm CDiv}}

\def\Supp{{\rm Supp}}

\def\cone{{\rm cone}}

\begin{document}

\title[Orbit spaces of small tori]
      {Orbit spaces of small tori}

\author[A.~A'Campo-Neuen]{Annette A'Campo-Neuen} 
\address{Fachbereich Mathematik, Johannes-Gutenberg-Universit\"at Mainz,
  55099~Mainz, Germany}
\email{acampo@enriques.mathematik.uni-mainz.de}
%
\author[J.~Hausen]{J\"urgen Hausen} 
\address{Fachbereich Mathematik und Statistik, Universit\"at Konstanz,
  78457 Konstanz, Germany}
\email{Juergen.Hausen@uni-konstanz.de}





\begin{abstract}
Consider an algebraic torus of small dimension acting on 
an open subset of $\CC^n$, or more generally on a quasiaffine
variety such that a separated orbit space exists. We discuss
under which conditions this orbit space is quasiprojective.
One of our counterexamples provides a toric variety with enough 
effective invariant Cartier divisors that is not embeddable
into a smooth toric variety.
\end{abstract}

\maketitle

\section*{Introduction}

Suppose a reductive group $G$ acts on a quasiprojective variety $X$
with a geometric quotient $X \to X / G$. The question when the
quotient space $X/G$ is again quasiprojective is studied by several 
authors, see e.g. \cite{BBSw2} and the classical counterexample 
\cite{BBSw1}. In the present note, we complement known partial 
answers for actions of low-dimensional tori on quasiaffine varieties.

A case of particular interest are diagonal actions of tori $T$ on 
$\CC^{n}$. Any maximal open subset $X \subset \CC^n$ admitting 
a geometric quotient $X \to X/T$ is in fact invariant under the
big torus $\TT^n := (\CC^*)^n$ and the orbit space $X/T$ is a 
simplicial toric variety with at most $n- \dim(T)$ invariant
prime divisors. We obtain: 

\begin{introthm}
Let $T$ be a torus acting diagonally on $\CC^{n}$, and let $X \subset
\CC^{n}$ be a $\TT^{n}$-invariant subset admitting a geometric
quotient $X \to X/T$. The following table indicates when $X/T$ 
is quasiprojective:

\bigskip

\begin{center}

\begin{tabular}{l||p{2.9cm}|p{2.9cm}|p{2.9cm}}
 \hfill & 
$X/T$ complete, \par $T$ acts freely & 
$X/T$ complete & 
$X/T$ arbitrary   \\
\hline \hline
$\dim(T) = 1$ & 
$+$ & 
$+$ & 
$+$, \cite[Example 5A]{BBSw2}   \\
\hline
$\dim(T) = 2$ & 
$+$, \cite[Theorem~2]{Kl} & 
$+$, Proposition~\ref{dimT2XCncomplete} & 
$-$, \cite[Example 5B]{BBSw2}  \\
\hline
$\dim(T) = 3$ & 
$+$, \cite[Theorem~1]{KlSt} & 
$-$, Proposition~\ref{dimT3XCncomplete} & 
$-$  \\
\hline
$\dim(T) = 4$ & 
$-$, \cite[Prop.~9.4]{Oda} & 
$-$  & 
$-$  
\end{tabular}

\end{center}

\bigskip

\noindent
Here the top row and the left column specify the assumptions on the
action, ``$+$'' stands for ``the quotient space $X/T$ is
quasiprojective'', and ``$-$'' means that $X/T$ is not necessarily
quasiprojective.
\end{introthm}

Let us turn to torus actions on arbitrary quasiaffine varieties 
$X$. If the torus $T$ is onedimensional and $X/T$ is complete, 
then the action defines a w.l.o.g. positive grading of the algebra
$\mathcal{O}(X)$ and the quotient space is nothing but 
${\rm Proj}(\mathcal{O}(X))$. As soon as we drop either of these two
assumptions, the orbit space is in general no longer quasiprojective,
see Propositions~\ref{dimT1Xquasiaffine} and~\ref{dimT2Xarbitcomplete}:

\goodbreak

\begin{introthm}
\begin{enumerate}
\item There is an action of a 2-dimensional torus $T$ on a
   5-dimen\-sio\-nal quasiaffine toric variety $X$ with a complete
   nonprojective orbit space $X/T$.
\item There is an action of a 1-dimensional torus $T$ on a
   5-dimensional quasiaffine toric variety $X$ with a
   nonquasiprojective orbit space $X/T$. 
 \end{enumerate}
\end{introthm}

The orbit space $X/T$ we construct to prove part ii) 
serves also as a subtle example in the context of
embeddings into toric varieties, compare~\cite{Wl} and~\cite{ha}: We
ask, which varieties can be embedded into smooth toric varieties. Note
that embeddability into a smooth toric variety requires existence of
``many'' Cartier divisors.

T.~Kajiwara~\cite{Ka} says that a toric variety has {\it enough
effective invariant Cartier divisors}, if the complements of these
divisors provide an affine cover. By~\cite{ha}, such a toric variety
always admits an embedding into a smooth toric prevariety with affine
diagonal, even by means of a toric morphism. However, we have:

\begin{introthm}
There exists a toric with enough effective invariant Cartier divisors
that admits no embedding into a separated smooth toric variety.%
\end{introthm}

This note is organized as follows: In Sec\-tion~\ref{dim3} we consider
actions of 3-dimen\-sio\-nal tori, and Sec\-tion~2 is devoted to 
2-dimensional torus actions. Finally, in Sections~\ref{kdiv} 
and~\ref{nontwodivisorialexample}, we construct the example of 
Theorem~2~ii) and Theorem~3.

\section{Quotients of threedimensional torus actions}\label{dim3}

Throughout the whole note, we work over the field $\CC$ of complex
numbers. First we recall the definition of a geometric quotient.
Consider an algebraic torus action $T \times X \to X$ on a complex
algebraic variety $X$. 

A  {\it good quotient} for this $T$-action is an affine $T$-invariant
regular map $p \colon X \to Y$ onto a variety $Y$ such that the
canonical map $\mathcal{O}_{Y} \to p_{*}(\mathcal{O}_{X})^{T}$ is an
isomorphism.  A good quotient $p \colon X \to Y$ of the $T$-action is
called {\it geometric} if the fibres $p^{-1}(y)$, $y \in Y$, are
precisely the $T$-orbits of $X$.

Given such a geometric quotient $X \to Y$, the variety $Y$ is called
the {\it orbit space} and is denoted by $X/T$. Note that
by~\cite[Corollary 3]{Su} and \cite[Propositions~0.7 and~0.8]{Mu}, an
effective algebraic torus action admits a geometric quotient if
and only if it is proper.

In this section, we consider the standard action of $\TT^6 = (\CC^*)^6$
on $\CC^6$ and present a 3-dimensional subtorus $T \subset \TT^6$
 and an open $\TT^6$-invariant subset $X \subset \CC^n$
with a complete but nonprojective orbit space $X/T$.

We shall use some basic notions of the theory of toric varieties; 
standard references are~\cite{Fu} and \cite{Oda}. All we need is the
following variant of a well-known example~\cite[Proposition~9.4]{Oda}:

\begin{example}\label{Odasexample}
Let $e_{1}$, $e_{2}$ and $e_{3}$ denote the canonical basis vectors
of the lattice $\ZZ^{3}$. Consider the vectors
  $$\begin{array}{ll}
    v_{1} := (-1,0,0), & \qquad v_{4} := (0,1,1), \\
    v_{2} := (0,-1,0), & \qquad v_{5} := (1,0,1), \\
    v_{3} := (0,0,-1), & \qquad v_{6} := (1,1,0). \\
\end{array}$$
Let $\Delta$ be the fan in $\ZZ^3$ with eight maximal cones, 
namely $\sigma_{i,j,k}:=\cone(v_i,v_j,v_k)$, where the triple
$(i,j,k)$ runs through the list
$$ (1,4,6), (1,3,6), (3,5,6), (2,3,5), (2,4,5),
(1,2,4), (1,2,3), (4,5,6).
$$

\bigskip

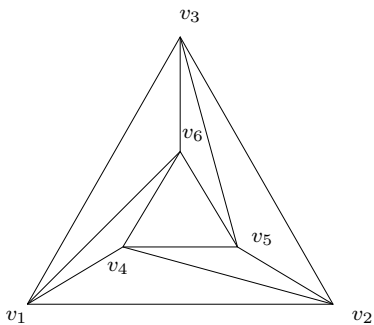
\begin{center}
  \begin{picture}(0,0)%
\includegraphics{odabeispielvar.pstex}%
\end{picture}%
\setlength{\unitlength}{0.00041700in}%
%
\begingroup\makeatletter\ifx\SetFigFont\undefined
\def\x#1#2#3#4#5#6#7\relax{\def\x{#1#2#3#4#5#6}}%
\expandafter\x\fmtname xxxxxx\relax \def\y{splain}%
\ifx\x\y   
\gdef\SetFigFont#1#2#3{%
  \ifnum #1<17\tiny\else \ifnum #1<20\small\else
  \ifnum #1<24\normalsize\else \ifnum #1<29\large\else
  \ifnum #1<34\Large\else \ifnum #1<41\LARGE\else
     \huge\fi\fi\fi\fi\fi\fi
  \csname #3\endcsname}%
\else
\gdef\SetFigFont#1#2#3{\begingroup
  \count@#1\relax \ifnum 25<\count@\count@25\fi
  \def\x{\endgroup\@setsize\SetFigFont{#2pt}}%
  \expandafter\x
    \csname \romannumeral\the\count@ pt\expandafter\endcsname
    \csname @\romannumeral\the\count@ pt\endcsname
  \csname #3\endcsname}%
\fi
\fi\endgroup
\begin{picture}(4824,3624)(4389,-3573)
\put(5431,-3136){\makebox(0,0)[lb]{\smash{\SetFigFont{8}{9.6}{rm}$v_4$}}}
\put(4141,-3751){\makebox(0,0)[lb]{\smash{\SetFigFont{8}{9.6}{rm}$v_1$}}}
\put(6346, 74){\makebox(0,0)[lb]{\smash{\SetFigFont{8}{9.6}{rm}$v_3$}}}
\put(6386,-1441){\makebox(0,0)[lb]{\smash{\SetFigFont{8}{9.6}{rm}$v_6$}}}
\put(7266,-2766){\makebox(0,0)[lb]{\smash{\SetFigFont{8}{9.6}{rm}$v_5$}}}
\put(8541,-3751){\makebox(0,0)[lb]{\smash{\SetFigFont{8}{9.6}{rm}$v_2$}}}
\end{picture}

\end{center}

\bigskip

\noindent%
Then the toric variety $X$ corresponding to the fan $\Delta$ is
simplicial and complete, but not projective.
\end{example}

\begin{proposition}\label{dimT3XCncomplete}
There exists a threedimensional subtorus $T \subset \TT^{6}$ and an
open $\TT^{6}$-invariant subset $X \subset \CC^{6}$ admitting a
geometric quotient with a nonprojective complete orbit space $X/T$.
\end{proposition}

\proof Using Cox's Construction \cite[Theorem~2.1]{Co}, we represent
the toric variety $X$ of Example~\ref{Odasexample} as a geometric
quotient of an open toric subvariety $\t{X} \subset \CC^{6}$: Define a
lattice homomorphism $Q\colon \ZZ^6\to\ZZ^3$ by $Q(e_i) := v_i$, and
for $\sigma \in \t{\Delta}$ consider
$$ \t{\sigma} := \cone(e_i; \; v_i \in \sigma). $$

These cones form a fan $\t{\Delta}$ in $\ZZ^{6}$. 
The associated toric variety $\t{X}$ is an open
toric subvariety of $\CC^6$ and the toric morphism $\t{X} \to X$
defined by $Q \colon \ZZ^{6} \to \ZZ^{3}$ is a geometric quotient for
the action of the subtorus $T \subset \TT^{6}$ defined by the
sublattice $\ker(Q) \subset \ZZ^{6}$.
\endproof

\section{Quotients for twodimensional torus actions}\label{dim2}

In this section we prove the statements on actions of twodimensionial
tori made in Theorems 1~and~2. The first result settles the case of 
actions on open subsets of $\CC^{n}$. As before, let 
$\TT^{n} := (\CC^{*})^{n}$, and endow $\CC^n$ with the standard 
$\TT^n$-action.

\begin{proposition}\label{dimT2XCncomplete}
Let $T \subset \TT^{n}$ be a subtorus of dimension two, and
suppose that $X \subset \CC^{n}$ is a $\TT^{n}$-invariant open
subset with a geometric quotient $X \to X/T$. If $X/T$ is complete, 
then $X/T$ is projective.  
\end{proposition}

Note that in the setting of this proposition, the orbit space $X/T$ is
a complete simplicial toric variety of dimension $n-2$ with at most $n$
invariant prime divisors. Thus the statement is an immediate consequence
of the following result:

\begin{proposition}\label{kleinschmidt}
Let $N$ be an $n$-dimensional lattice, and let $\Delta$ be a complete
simplicial fan in $N$ having at most $n+2$ rays. Then $\Delta$ is
strongly polytopal.
\end{proposition}

Recall that P. Kleinschmidt proved this for the case that $\Delta$ is a
regular fan, see \cite{Kl}.
As in \cite{Kl}, we shall use the following projectivity criterion in
terms of Gale transforms, first stated by Shephard and later employed
by Ewald, compare~\cite{Ew} and \cite{Sh}:

Suppose that $R$ denotes a list of  generators of the $d$ rays of a
complete fan $\Delta$ in $\ZZ^n$. A {\it linear Gale transform\/}
$\b{R}$ of $R$ consists of the columns of a matrix whose rows form
a basis for the linear relations of $R$. The {\it coface\/} of a cone
$\sigma\in\Delta$ generated by $v_{i_1},\dots,v_{i_r}$ is
$$\b{\sigma} := 
\cone(\b{R}\setminus\{\b{v_{i_1}},\dots,\b{v_{i_r}}\}) 
  \subset\RR^{d-n}.$$ 

\begin{lemma} 
The fan $\Delta$ is strongly polytopal if and only if
the intersection of the relative interiors of all $\b{\sigma}$,
$\sigma \in \Delta^{\max}$, is nonempty.
\end{lemma}

\medskip

\proof[Proof of Proposition~\ref{kleinschmidt}]
If $\Delta$ has $n+1$ rays, then the primitive vectors in the rays
generate a simplex and there is nothing to show. So let us  assume
that $\Delta$ has $n+2$ rays. First we reduce to the case that there
is a regular maximal cone $\sigma_{0} \in \Delta$: 

Fix any maximal cone $\sigma_0$ of $\Delta$ and let $N_0 \subset N$ be
the sublattice spanned by the primitive generators of $\sigma_0$. The
cones of $\Delta$ also form a fan $\Delta_{0}$ in $N_0$. Clearly, if
$\Delta_{0}$ is strongly polytopial, then so is $\Delta$. Thus,
replacing $\Delta$ with $\Delta_{0}$, we may assume $\sigma_{0} \in
\Delta$ is regular. 

Choosing the primitive generators of $\sigma_0$ as a basis, we may
assume that $N = \ZZ^{n}$ and that $\sigma_{0}$ is generated by the
canonical base vectors $e_{1}, \ldots, e_{n}$. We denote the primitive
vectors of the remaining two rays of $\Delta$ by $u$ and $v$.
By the combinatorial classification of spherical complexes, see
\cite{Ma} and \cite{Gr}, the list  
$$R:=(e_1,\dots,e_n,u,v)$$
can be partitioned into two complementary subsets $U$ and $V$ with
$u\in U$, $v\in V$ and $|U|,|V|\geq 2$ such that 
$$ \Delta^{\max}=\{\cone(R\setminus \{w,z\}); \; w\in U,z\in V\}\,.$$

After renumbering, we may assume that $U$ equals $\{e_1,\dots,e_r,u\}$
and $V$ equals $\{e_{r+1},\dots,e_n,v\}$ for some $1\leq r<n$. Then,
in addition to $\sigma_0$, we have the following three types of
maximal cones in $\Delta$:
$$ \begin{array}{ll}
\sigma_{i}:=\cone(u,e_k; \; k\ne i) 
   & \quad \text{for $1\leq i\leq r$}, \\
\sigma_{j}:=\cone(v,e_k; \; k\ne j) 
   & \quad \text{for $r < j  \leq n$}, \\ 
\sigma_{ij}=\cone(u,v,e_k; \; k\ne i,j) 
   & \quad \text{for $1 \leq i \leq r$ and $r < j \leq n$}.
\end{array}
$$

Note that the linear form $e_{i}^{*}$ separates the pair
$\sigma_{0}$, $\sigma_{i}$, i.e.~the two cones lie on different sides
of the hyperplane defined by $e_{i}^{*}$. 
Similarly, $e_{j}$ separates $\sigma_{0}$,
$\sigma_{j}$. Moreover, $\sigma_{i}$ and $\sigma_{ij}$ are separated
by $u_{i}e_{j}^{*}-u_{j}e_{i}^{*}$. Applying these linear forms to the
generators of the respective cones successively yields
$$ \begin{array}{ll}
u_{i} < 0  & \quad \text{for $1\leq i\leq r$}, \\
v_{j} < 0  & \quad \text{for $r < j  \leq n$}, \\ 
u_{i}v_{j} - u_{j}v_{i} > 0
   & \quad \text{for $1 \leq i \leq r$ and $r < j \leq n$}.
\end{array}
$$ 

We shall combine these inequalities with the above mentioned
projectivity criterion. Consider the following linear Gale transform
of $R$:
$$\b{R}=((u_1,v_1),\dots,(u_n,v_n),(-1,0),(0,-1)).$$ 
Then the cofaces associated to the maximal cones of $\Delta$ are 
given by 
$$
\begin{array}{ll}
\b{\sigma_{0}} := \cone((-1,0),(0,-1)), & \\
\b{\sigma_{i}} := \cone((u_i,v_i),(0,-1)) 
   & \quad \text{for $1\leq i\leq r$}, \\
\b{\sigma_{j}} := \cone((-1,0),(u_j,v_j)) & 
   \quad \text{for $j < r \leq n$}, \\ 
\b{\sigma_{ij}} := \cone((u_i,v_i),(u_j,v_j)) &
   \quad \text{for $1\leq i\leq r$ and $r<j\leq n$}.
\end{array}$$

The above inequalities imply that all these cones are twodimensional,
and by Shephard's criterion we only have to show their intersection is
again twodimensional. There are $1 \le i_{0} \le r$ and $r < j_{0} \le n$
with 
$$\b{\sigma_{i_0}} = \bigcap_{i\leq r}\b{\sigma_{i}}, 
\qquad \b{\sigma_{j_0}} = \bigcap_{j> r}\b{\sigma_{j}}\, $$
and the intersection over all cones $\b{\sigma_{i,j}}$ equals
$\b{\sigma_{i_{0},j_{0}}}$. Consequently we obtain for the
intersection of all cofaces of maximal cones:
$$ \bigcap_{\sigma\in \Delta^{\max}} \b{\sigma} 
= \cone((-1,0),(0,-1)) \cap \b{\sigma_{i_0,j_0}}. $$
Using  the above inequalities for $i_0$ and $j_0$, 
we conclude that this cone is in
fact of dimension two. \qed

\medskip

Now we turn to an arbitrary quasiaffine variety $X$. We show by means
of an example that the orbit space $X/T$  in general may not be
quasiprojective:

\begin{proposition}\label{dimT2Xarbitcomplete}
There exists a five dimensional quasiaffine toric variety $X$ with big
torus $T_{X}$ and a twodimensional subtorus $T \subset T_{X}$ acting
with geometric quotient on $X$ such that $X/T$ is a complete
nonprojective variety. 
\end{proposition}

\proof We will represent the threedimensional toric variety $X$
introduced in Example~\ref{Odasexample} as an orbit space of the
action of a twodimensional torus on a quasiaffine toric variety
$\rq{X}$. Consider the following vectors in $\ZZ^5$:
$$\begin{array}{ll}
    w_{1} := (1,0,0,0,0), & \qquad w_{4} := (0,0,1,0,0), \\
    w_{2} := (0,1,0,0,0), & \qquad w_{5} := (0,0,0,1,0), \\
    w_{3} := (-1,0,0,-1,2), & \qquad w_{6} := (0,-2,-1,1,1). \\
\end{array}$$
Let $v_{i} \in \ZZ^{3}$ as in Example~\ref{Odasexample}. Then there is
a lattice homomorphism $P\colon \ZZ^5\to\ZZ^3$ with $P(w_i)=v_i$ for
all $i$, namely the homomorphism defined by the matrix
$$
\left[ 
{\begin{array}{rrrrr}
-1 & 0 & 0 & 1 & 0 \\
0 & -1 & 1 & 0 & 0 \\
0 & 0 & 1 & 1 & 0
\end{array}}
 \right] \,.
$$

One directly checks that the vectors $w_i$ generate a strictly convex
cone $\sigma_0$ in $\ZZ^5$. Moreover, for every maximal cone
$\sigma\in\Delta$, one obtains a face $\rq{\sigma}$ of $\sigma_{0}$ by
setting
$$\rq{\sigma} := \cone(w_i; \; v_i\in\sigma). $$

The fan $\rq{\Delta}$ generated by these $\rq{\sigma}$
corresponds to a quasiaffine toric variety $\rq{X}$. Since
$\dim\rq{\sigma}=\dim\sigma$ for all $\sigma\in\Delta$, the
toric morphism $p\colon \rq{X}\to X$ defined by the lattice
homomorphism $P \colon \ZZ^{5} \to \ZZ^{3}$ is the desired   geometric
quotient, use e.g., \cite[Theorem~5.1]{Hm}. \qed

\section{When is a toric variety $k$-divisorial?}\label{kdiv}

In this section we give a criterion for $k$-divisoriality, a notion
that comes up naturally in the context of toric embeddings.
The criterion will be used in the following section for the
construction of a toric variety with enough effective invariant
Cartier divisors that cannot be embedded into
a smooth toric variety.

As in \cite{ha}, we call an irreducible variety $X$ {\it
$k$-divisorial} if any $k$-points $x_{1}, \ldots, x_{k}$ admit a
common affine neighbourhood of the form $X \setminus \Supp(D)$ with
an effective Cartier divisor $D$ on $X$. For $k=1$ this gives back the
usual notion of a divisorial variety, see e.g.~\cite{Bo}, \cite{SGA}.

In the whole section we assume that our toric variety
$X$ is non degenerate, i.e., there exists no toric decomposition $X
\cong X' \times \KK^{*}$. Our criterion reads as follows:

\begin{proposition}\label{kdiviscriterion2}
A toric variety $X$ is $k$-divisorial if and only if for any $k$
closed orbits $B_{1}, \ldots, B_{k} \subset X$ of the big torus $T_{X}
\subset X$ there exist $T_{X}$-invariant effective Cartier divisors
$D_{1}, \ldots, D_{k}$ on $X$ such that 
\begin{enumerate}
\item any two $D_{i}$, $D_{j}$ are linearly equivalent to each other,
\item $X \setminus \Supp(D_{i})$ is an affine neighbourhood of
  $B_{i}$.
\end{enumerate}
\end{proposition}

Note that for $k=1$, this is just~\cite[Proposition~1.2]{acha2} for
toric varieties. Moreover, if we choose the number $k$ in this
proposition to be the number of all closed $T_{X}$-orbits, then we
obtain the following well known result, compare \cite{Su}, Lemma~8:

\begin{corollary}
A toric variety $X$ with at most $k$ closed $T_{X}$-orbits is 
quasiprojective if and only if it is $k$-divisorial. \qed
\end{corollary}

The proof of Proposition~\ref{kdiviscriterion2} relies on a more
explicit formulation of the result which we state below in 
Proposition~\ref{kdiviscriterion}. For this, we 
use Cox's construction \cite[Theorem~2.1]{Co} to obtain an open toric
subset $\t{X}$ of some $\CC^{n}$ with complement $\CC^{n} \setminus
\t{X}$ of dimension at most $n-2$ and a closed subgroup $H \subset 
\TT^{n}$ such that there is good quotient
$$ q \colon \t{X} \to \t{X} \quot H = X. $$
Note that this quotient map is a toric morphism. Suppose now that $X$
is divisorial. Then Kajiwara's construction \cite[Theorem 1.9]{Ka}
gives rise to a commutative diagram of toric morphisms 
$$ \xymatrix{
{\t{X}} \ar[rr]_{q_{1}}^{\quot H_{1}} \ar[dr]^{q}_{\quot H} & &
{\rq{X}} \ar[dl]_{q_{2}}^{/ H_{2}} \\
& X & }$$
where $\rq{X}$ is a quasiaffine toric variety, $H_{1}, H_{2}$ are
closed subgroups of the respective big tori $T_{\t{X}}$ and
$T_{\rq{X}}$ and the map $q_{2} \colon \rq{X}
\to X$ is a geometric quotient. In our criterion, we use the notion
of a {\it distinguished point}, see \cite[p.~28]{Fu}. The distinguished
points of $\CC^{n}$ are just the points having only coordinates $0$ or
$1$.

\begin{proposition}\label{kdiviscriterion}
The toric variety $X$ is $k$-divisorial if and only if for any $k$
distinguished points $\t{x}_{1}, \ldots, \t{x}_{k} \in \t{X}$ with
$\TT^{n} \mal x_{i}$ closed in $\t{X}$ there exist monomials $f_{1}, \ldots,
f_{k} \in \CC[z_{1}, \ldots, z_{n}]$ such that
\begin{enumerate}
\item the polynomial $f := f_{1} + \ldots + f_{r}$ is $H$-homogeneous
  and $H_{1}$-invariant,
\item $f_{i}(\t{x}_{i}) = 1$ and $f_{i}(\t{x}) = 0$ for every $\t{x} \in
  \b{\TT^{n} \mal \t{x}_{i}} \setminus \TT^{n} \mal \t{x}_{i} \subset
  \CC^{n}$.
\end{enumerate}
\end{proposition}

\proof Assume first that $X$ is $k$-divisorial and let $\t{x}_{1},
\ldots, \t{x}_{k} \in \t{X}$ be distinguished points with closed
$\TT^{n}$-orbit in $\t{X}$. Then there exists an effective Cartier
divisor $D \in \CDiv(X)$ such that $X \setminus \Supp(D)$ is affine
and contains the points $q(\t{x}_{1}), \ldots, q(\t{x}_{k})$.

We claim that the pullback $q_{2}^{*}(D) \in \CDiv(\rq{X})$ is
principal. To see this note first that $D$ is of the form $E + \div(h)$ 
with some $T_{X}$-invariant Cartier divisor $E \in \CDiv(X)$ and some
function $h \in \CC(X)$, see e.g. \cite[p. 63]{Fu}. Thus we obtain
$$ q_{2}^{*}(D) = q_{2}^{*}(E) + \div(q_{2}^{*}(h)). $$
By \cite[Proposition~2.6]{acha1}, we have $q_{2}^{*}(E) = \div(\rq{g})$ 
with some character function $\rq{g}$ of the big torus $T_{\rq{X}}
\subset \rq{X}$. Setting $\rq{h} := \rq{g}q_{2}^{*}(h)$, we have
$q_{2}^{*}(D) = \div(\rq{h})$. Note that $\rq{h}$ is $H_{2}$-homogeneous
and, since $q_{2}^{*}(D)$ is effective, $\rq{h}$ is a regular function
on $\rq{X}$. In particular, our claim is verified.

To proceed, consider $\t{h} := q_{1}^{*}(\rq{h})$. This is a regular
function on $\t{X}$ and hence it is a polynomial. Moreover, $\t{h}$ is
$H_{1}$-invariant and $H$-homogeneous. Finally, on $\t{X}$, we
have $\div(\t{h}) = q^{*}(D)$. Since $q^{-1}(X \setminus
\Supp(D))$ is affine and $\CC^{n} \setminus \t{X}$ is small, this
implies $\t{X}_{\t{h}} = \CC^{n}_{\t{h}}$. 

Now, consider one of the distinguished points $\t{x}_{i}$. By
construction, we have $\t{h}(\t{x}_{i}) \ne 0$. Consequently, 
there appears
a monomial $f_{i}$ in $\t{h}$ with $f_{i}(\t{x}_{i}) \ne 0$. Surely, also
$f_{i}$ is $H_{1}$-invariant and it is $H$-homogeneous with respect to
the same weight as $\t{h}$. Let
$$ f_{i}(z_{1}, \ldots, z_{n}) = z_{1}^{u_{1}} \ldots z_{n}^{u_{n}},
\qquad  \t{x}_{i} = (\t{x}_{i1}, \ldots, \t{x}_{in}).$$

Then clearly $\t{x}_{ij} = 0$ implies $u_{j} = 0$. We claim that also
the converse is true. Indeed, suppose that $u_{j} = 0$ but $x_{ij}
\neq 0$ holds for some $j$. Then, replacing in $x_{i}$ the coordinate
$x_{ij}$ with zero yields a point
$$x_{i}' \in \CC^{n} \setminus \t{X}.$$
with $f_{i}(x_{i}') \ne 0$. This implies that the restriction of
$\t{h}$ to the orbit $\TT^{n} \mal x_{i}'$ is not the zero
function. But this contradicts the fact that $\t{X}_{\t{h}}$
equals $\CC^{n}_{\t{h}}$. So our claim is verified.

Now, take for each distinguished point $\t{x}_{i}$ a monomial $f_{i}$ as
above. Then these monomials fullfill the desired conditions, and one
implication of the proposition is proved.

For the reverse direction, suppose that conditions i) and ii) of the
assertion hold. First we consider $k$ distinguished points $x_{1},
\ldots, x_{k} \in X$ with $T_{X} \mal x_{i}$ closed in $X$. Choose
distinguished point $\t{x}_{1}, \ldots, \t{x}_{k} \in \t{X}$ such that
$x_{i} = q(\t{x}_{i})$ holds and $\TT^{n} \mal \t{x}_{i}$ is closed in
$\t{X}$. 

Let $f_{1}, \ldots, f_{k}$ be polynomials satisfying conditions i) and
ii), and let $f := f_{1} + \ldots + f_{k}$. Since $f$ is
$H_{1}$-invariant, it is of the form $f = q_{2}^{*}(\rq{h})$ for some $\rq{h}
\in \mathcal{O}(\rq{X})$. It follows as in \cite[Proof of
Proposition~1.3]{acha2} that there is an effective Cartier divisor $D$
on $X$ with 
$$\Supp(D) = q_{2}(\Supp(\div(\rq{h}))) = q(\Supp(\div(f))).$$
By the properties of the $f_{i}$, we have $f(\t{x}_{i}) =1$. Hence the
above equation yields $x_{1}, \ldots, x_{k} \in X \setminus
\Supp(D)$. This settles the case of $k$ distinguished points $x_{1},
\ldots, x_{k}$ with closed $T_{X}$-orbits.

If $x_{1}', \ldots, x_{k}' \in X$ are arbitrary, then we
can choose distinguished points $x_{i}$ in the closure of $T_{X} \mal
x'_{i}$ with $T_{X} \mal x_{i}$ closed in $X$. The above consideration
provides an effective Cartier divisor $D$ on $X$ such that $U := X \setminus
\Supp(D)$ is an affine neighbourhood of $x_{1}, \ldots, x_{k}$. Let $t
\in T_{X}$ with $t \mal x_{i}' \in U$. Then $t^{-1} \mal D$ is the
desired Cartier divisor. \qed 

\medskip

\proof[Proof of Proposition~\ref{kdiviscriterion2}] 
In the setting of Proposition~\ref{kdiviscriterion}, the monomials
$f_{i}$ correspond to Cartier divisors $D_{i}$ on $X$ and vice versa,
see e.g. \cite[Proposition~2.6]{acha1}. The properties of the $f_{i}$
translate directly to the desired properties of the $D_{i}$. Here
the fact that the $D_{i}$ are pairwise linearly equivalent,
corresponds to the  fact that all $f_{i}$ are $H$-homogeneous with
respect to the same weight. \qed

\section{A divisorial toric variety that is not 2-divisorial}%
\label{nontwodivisorialexample}

In this section we present an example of a fourdimensional divisorial
toric variety $X$ that is not 2-divisorial. In order to define the fan
of $X$, consider the following lattice vectors in $\ZZ^4$:
$$ \begin{array}{lll}
v_{1} := (1,0,0,0), & v_{2} := (0,-2,1,0), & v_{3} := (0,-1,1,1), \\
v_{4} := (0,0,0,1), & v_{5} := (0,1,0,2),  & v_{6} := (-1,-1,-1,2), \\
v_{7} := (1,-1,0,-1), & v_{8} := (1,1,-1,0), & v_{9} := (0,0,-2,1).
\end{array}
$$
The toric variety $X$ is defined by the fan $\Delta$ in $\ZZ^4$ having the
following five cones as its maximal cones:
$$\begin{array}{lll}
\sigma_1 := \cone(v_2,v_3,v_4,v_5,v_6), & 
\sigma_2 := \cone(v_1,v_2,v_3,v_7), &
\sigma_3 := \cone(v_4,v_5,v_8), \\
\sigma_4 := \cone(v_1,v_5,v_8), &
\sigma_5 := \cone(v_1,v_7,v_9).
\end{array}
$$

\begin{proposition}\label{nontwodivisorialexampleprop}
The toric variety $X$ defined by the fan $\Delta$ is divisorial but not
2-divisorial.
\end{proposition}

\proof In order to apply Proposition~\ref{kdiviscriterion}, we have to
determine the quotient presentations of $X$ due to 
Cox and Kajiwara. To obtain Cox's construction,
consider the fan $\t{\Delta}$ in $\ZZ^9$ generated by the cones
$$ \t{\sigma}_i := \cone(e_j; \; v_j \in \sigma_i), \qquad i=1,\ldots,
5.$$
The associated toric variety $\t{X}$ is an open toric subvariety of
$\CC^{9}$ with $7$-dimensional complement. Moreover, the lattice
homomorphism $Q\colon \ZZ^9\to\ZZ^4$ sending the canonical base vector
$e_i$ to $v_i$, induces Cox's quotient presentation $q \colon \t{X}\to
X$.

To obtain Kajiwara's quotient presentation, we have to determine the
invariant Cartier divisors of $X$. Let $D_{i}$ denote the
invariant Weil divisor on $X$ corresponding to the ray through
$v_{i}$. An explicit calculation shows that
$$ D_{6}, \ D_{8}, \ D_{9}, \ D_{1}+D_{7}, \ D_{2}+D_{3}, \
D_{4}+D_{5}, \ D_{2}+D_{4}+D_{7}  $$
form a basis for the group of invariant Cartier divisors of $X$. Thus
we have to consider the lattice homomorphism $Q_{1} \colon \ZZ^{9} \to
\ZZ^{7}$ given by the matrix
$$ \left[ 
\begin{array}{rrrrrrrrr}
0 & 0 & 0 & 0 & 0 & 1 & 0 & 0 & 0 \\
0 & 0 & 0 & 0 & 0 & 0 & 0 & 1 & 0 \\
0 & 0 & 0 & 0 & 0 & 0 & 0 & 0 & 1 \\
1 & 0 & 0 & 0 & 0 & 0 & 1 & 0 & 0 \\
0 & 1 & 1 & 0 & 0 & 0 & 0 & 0 & 0 \\
0 & 0 & 0 & 1 & 1 & 0 & 0 & 0 & 0 \\
0 & 1 & 0 & 1 & 0 & 0 & 1 & 0 & 0
\end{array}
\right]
$$
One directly checks that the vectors $Q_{1}(e_{1}), \ldots,
Q_{1}(e_{9})$ generate the extremal rays of a strictly convex cone
$\rq{\sigma} \subset \QQ^{7}$ and that the cones
$$ \rq{\sigma}_{i} := Q_{1}(\t{\sigma_{i}}), \qquad i = 1, \ldots, 5$$
are faces of $\rq{\sigma}$. Consequently, the fan $\rq{\Delta}$ in $\ZZ^{7}$
generated by $\rq{\sigma}_{1}, \ldots, \rq{\sigma}_{5}$ defines a
quasiaffine toric variety $\rq{X}$. Moreover, we have a commutative
diagram of toric morphisms
$$ \xymatrix{
{\t{X}} \ar[rr]_{q_{1}}^{\quot H_{1}} \ar[dr]^{q}_{\quot H} & &
{\rq{X}} \ar[dl]_{q_{2}}^{/ H_{2}} \\
& X & }$$
where $q_{1} \colon \t{X} \to \rq{X}$ arises from the lattice
homomorphism $Q_{1} \colon \ZZ^{9} \to \ZZ^{7}$. In particular, $X$ is
divisorial and the toric morphism $q_{2} \colon \rq{X} \to X$ is
Kajiwara's quotient presentation.


Now assume that $X$ were even $2$-divisorial. The
distinguished points $\t{x}_{3},\t{x}_{5} \in \t{X}$ corresponding to
the maximal cones $\t\sigma_3$ and $\t\sigma_5$ of $\t\Delta$ are
given by 
$$\t{x}_3 = (1,1,1,0,0,1,1,0,1), \qquad 
\t{x}_5 = (0,1,1,1,1,1,0,1,0). $$
Let $H \subset \TT^{9}$ and $H_{1} \subset \TT^{9}$ denote the subtori
corresponding to $\ker(Q)$ and $\ker(Q_{1})$
respectively. Proposition~\ref{kdiviscriterion} yields monomials of
the form
$$ f_3(z_{1}, \ldots, z_{9}) =
z_1^{a_1}z_2^{a_2}z_3^{a_3}z_6^{a_6}z_7^{a_7}z_9^{a_9}, 
\qquad
f_5(z_{1}, \ldots, z_{9}) =
z_2^{b_2}z_3^{b_3}z_4^{b_4}z_5^{b_5}z_6^{b_6}z_8^{b_8} $$
where $a_i,b_i>0$, such that $f_3$, $f_5$ are invariant with respect
to $H_{1}$ and both are $H$-homogeneous with respect to the same
weight. The kernel of $Q_1$ is generated by the vectors
$$ (-1,-1,1,0,0,0,1,0,0) \qquad (-1,0,0,-1,1,0,1,0,0). $$
Thus, $H_{1}$-invariance of the monomials $f_{3}$ and $f_{5}$ implies 
$$a_1=a_7, \quad a_2=a_3, \quad b_4=b_5,  \quad b_2=b_3.$$
Now consider the one parameter subgroup $\CC^{*} \to H$ corresponding
to the following lattice vector
$$(3,-1,1,-3,0,0,-1,-2,1) \in \ker(Q).$$
Since both $f_3$ and $f_5$ have the same weight with respect to $H$,
they also have the same weight with respect to the above one
parameter subgroup, and we obtain the relation
$$ 2 a_1+a_9=-3b_4-2b_8. $$
This contradicts the assumption that all the exponents $a_{i}$ and
$b_{i}$ are positive. Consequently, $X$ cannot be $2$-divisorial. \qed

\medskip

In \cite{ha}, 1- and 2-divisoriality are characterized in terms of
toric embeddings. 
Using~\cite[Theorems~3.1 and~3.2]{ha}, we obtain:

\begin{theorem}
The toric variety $X$ of Proposition~\ref{nontwodivisorialexampleprop}
admits a closed toric embedding into a smooth
toric prevariety with affine diagonal but it cannot be embedded into
a smooth toric variety.
\end{theorem}
 
We turn back to quasiprojectivity of orbit spaces. Our example serves
also to show that for a $\CC^{*}$-action on a quasiaffine variety with
geometric quotient, the resulting orbit space  in general need not be
quasiprojective: 

\begin{proposition}\label{dimT1Xquasiaffine}
There exists a five-dimensional affine toric variety $X$ with 
$\CC^{*}$-action and a toric open subset $U \subset X$ admitting a
geometric quotient $U \to U / \CC^{*}$ such that $U / \CC^{*}$ is not
$2$-divisorial.
\end{proposition}

\proof The affine toric variety $X$ in question arises from the
five-dimensional cone $\tau$ in $\ZZ^5$ generated by the following
$9$ lattice vectors:
$$\begin{array}{lll}
w_1:=(0,0,0,1,0), &
w_2:=(0,0,1,0,-1), &
w_3:=(0,0,1,0,0), \\
w_4:=(0,1,0,0,0), &
w_5:=(0,1,0,0,1), &
w_6:=(1,0,1,-1,0),\\
w_7:=(0,0,0,1,-1), &
w_8:=(0,1,-1,1,0), &
w_9:=(1,0,0,0,0).
\end{array}$$
We consider the action of the one parameter subgroup $\CC^{*} \to
\TT^{5}$ on $X$ corresponding to the lattice vector
$$ w := (1,-5,2,0,2). $$
The open toric subset $U \subset X$ is given by the fan $\Sigma$ in
$\ZZ^5$ with the following maximal cones:
$$\begin{array}{ll}
\tau_1 := \cone(w_2,w_3,w_4,w_5,w_6), &
\tau_2 := \cone(w_1,w_2,w_3,w_7), \\
\tau_3 := \cone(w_4,w_5,w_8), &
\tau_4 := \cone(w_1,w_5,w_8), \\
\tau_5 := \cone(w_1,w_7,w_9).
\end{array}
$$
Note that $\tau_{1}, \ldots, \tau_{5}$ are in fact faces
of the cone $\tau$. Let $P \colon \ZZ^5 \to \ZZ^4$ be the projection
defined by the matrix
$$\left[ 
{\begin{array}{rrrrr}
0 & 0 & 0 & 1 & 0 \\
0 & 0 & -1 & 0 & 1 \\
-2 & 0 & 1 & 0 & 0 \\
1 & 1 & 1 & 0 & 1
\end{array}}
\right]
$$
The kernel of $P$ is generated by the lattice vector $w \in
\ZZ^{5}$. Projecting the cones of $\Sigma$ via $P$, we
obtain just the fan $\Delta$ of the toric variety
presented in Proposition~\ref{nontwodivisorialexampleprop}.

Since $P$ is injective on all the cones $\tau_i$, it defines in fact a
geometric quotient $U \to U / \CC^{*}$ for the $\CC^*$-action on
$U$ corresponding to $w \in \ZZ^{5}$. As we have seen in
Proposition~\ref{nontwodivisorialexampleprop}, the quotient variety
$U/\CC^*$ is not $2$-divisorial. \qed

\bigskip

\noindent{\bf Acknowledgement.}\enspace
In the search for explicit examples we used the Maple-Package 
{\tt convex} developed by M.~Franz. The package contains several
basic functions on convex geometry. We would like to thank M.~Franz
for making his package available under
{\tt http://www.mathe.uni-konstanz.de/\~{ }franz/convex/}.

\bibliography{}

\end{document}